\documentstyle[12pt]{article}
\begin{document}
\begin{center}

{\bf SUPPORT OF BORELIAN MEASURES }

\vspace{3mm}

  {\bf IN SEPARABLE BANACH SPACES.}\\

\vspace{3mm}

{\sc Ostrovsky E.}\\

\vspace{3mm}

{\it Department of Mathematics and Statistics, Bar-Ilan University,
59200, Ramat Gan, Israel.}\\
e-mail: \ galo@list.ru \\

\vspace{5mm}

{\sc ABSTRACT}

\end{center}

\vspace{3mm}

 We prove in this article that every Borelian measure, for example, the
distribution of a random variable, in separable Banach space has a support
which is compact embedded Banach subspace; and prove that if the norm of the random variable belongs to some exponential Orlicz space, then the new subspace
can be choose such that the norm of this variable in the new space also belongs
to other exponential Orlicz space.

\vspace{3mm}

 {\it Key words:}  Banach and Orlicz space, Borelian measure, support, compact
embedded Banach subspace,  norm, rearrangement invariant space, metric entropy,
distance, moment, Grand Lebesgue Space, martingale. \\

{\it Mathematics Subject Classification (2000):} primary 60G17; \ secondary
 60E07; 60G70.\\

\vspace{3mm}

{\bf 1. Introduction. Notations. Statement of problem.}

\vspace{3mm}

 Let \\
{\bf A ) } $ \ X $ be Separable Banach Space (SBS) with a norm $ |\cdot|X, $
briefly:  $ X = (X, |\cdot|X) \in \ SBS, $  equipped with Borelian sigma \ -
\ algebra $ B; $ \\

{\bf B )} $ (\Omega, F, {\bf P}) $ be probability triple with expectation
$ {\bf E}; $ \\

{\bf C  ) } $ \ \xi $ be a random variable (r.v.) defined on the our triple
with the values in the space $ X:$

$$
{\bf P}(\xi \notin X) = 0 \eqno(1)
$$

with  Borelian distribution  $ \mu_{\xi}(\cdot): $

$$
\mu_{\xi}(A) = {\bf P}(\xi(\cdot) \in A), \ A \in B;
$$

{\bf D ) } $ \Phi = \Phi(u), \ u \in R $ be an Orlicz function, i.e. convex, even, continuous, twice continuous differentiable in the domain
$ u \in [2,\infty), $  strictly decreasing in the right side half line, and
such that

$$
\Phi(0) = 0; \ \lim_{u \to \infty} \Phi(u) = \infty.
$$

 We will denote the Orlicz space over $ (\Omega, F, {\bf P}) $ with Orlicz
function $ \Phi(\cdot) $ as $ Or(\Phi) $ and the correspondent Orlicz norm as
$ \ || \cdot ||Or(\Phi): $ for all {\it real } valued r.v. $ \eta $

$$
|| \eta ||Or(\Phi) = \inf \{\lambda, \lambda > 0, \ {\bf E}
\Phi(|\eta|/\lambda) \le 1  \}.
$$

 By definition, $ \eta \in Or(\Phi) $ iff $ ||\eta||Or(\Phi) < \infty. $ \par
It is proved ([1], [2], [5], [6]) that  if the condition (1) is
satisfied,  then there exists a compact embedded SBS in the space $ X $
 subspace $ Y = (Y, |\cdot|Y): \ Y \subset X $  such that the distribution of $ \xi $ is concentrated on the subspace $ Y: $

$$
{\bf P}(\xi \notin Y) = 0. \eqno(2)
$$

 Recall that the SBS subspace $ Y $ of the space $ X: \ Y \subset X $  is
called  compact embedded, if the unit ball of the space $ Y $ is precompact
set in the space $ X, $ i.e. that the closure in the $ X $ norm of unit
ball of the space $ Y $ is compact set in the space $ X. $ \par
 Moreover, if for some Orlicz function $ \Phi $ satisfying the so-called
$ \Delta_2 $  condition:

$$
\overline{\lim}_{u \to \infty} \frac{\Phi(2u)}{\Phi(u)} < \infty \eqno(3)
$$

the norm of r.v. $ \xi, $ belongs to Orlicz space $ Or(\Phi): $
$$
|| \ |\xi|X \ ||Or(\Phi) < \infty, \eqno(4)
$$

we can choose the subspace $ Y = (Y, |\cdot|Y) $  such that the {\it new} norm
of a r.v. $ \xi $ in the space $ Y, $  i.e. the variable  $ |\xi|Y $ belongs to
the space $ Or(\Phi): $

$$
|| \ |\xi|Y \ ||Or(\Phi) < \infty.
$$

{\bf The aim of this report is investigate the case when the function $ \Phi $
does not satisfies the $ \Delta_2 $ condition (3). } \par

 Recall in the end oh this section then the Orlicz function $ \Psi $ is called
weaker than $ \Phi, $ if  for all positive constant $ v; \ v = const > 0 $

$$
\lim_{u \to \infty} \frac{\Psi(v \ u )}{\Phi(u)} = 0.   \eqno(5)
$$

Notation: $ \Psi << \Phi. $ \par

\vspace{3mm}

{\bf 2. Main result.} \par

{\bf Theorem 1.} Let $ \xi $ be a r.v. with Borelian distribution
$ \mu_{\xi}(\cdot) $ in the SBS $ (X, | \cdot|X |) $ satisfies the condition (1). \par
 Let also $ \Psi(\cdot) $  be other Orlicz function weaker than $ \Phi(\cdot):
\Psi << \Phi. $ \par

 Then there exists a compact embedded in the space $ X $ the SBS space $ Y =
(Y, |\cdot|Y), $  which depended on the $ \mu_{\xi} $ and on the function
$ \Psi, $ such that the variable $ \xi $ belongs to the space $ Y $  a.e. and moreover

$$
|| \ |\xi|Y \ ||Or(\Psi) < \infty. \eqno(6)
$$

\vspace{3mm}

{\bf Proof.} \par
{\bf 1.} Note first of all that we can conclude that the space $  X $ coincides
with the SBS space $ C[0,1] $ of all continuous real valued function defined on
the closed interval [0,1], as long as the last space is universal in the class
of all SBS. \par

 We denote the norm in the space $ C[0,1] $  as $ |f|_{\infty}: $

$$
|f|_{\infty} = \sup_{t \in [0,1]} |f(t)|.
$$

{\bf 2.} In detail, the r.v. $ \xi $ can be identified with a continuous a.e. random process $  \xi = \xi(t), \ t \in [0,1]. $ \par
 Since the variable $ |\xi|_{\infty} $ belongs to the Orlicz space $ Or(\Phi), $
we can and will assume that

$$
{\bf E} \ \Phi \left( |\xi|_{\infty} \right) \le 1.
$$

 Let us define as usually for arbitrary continuous on the set [0,1] function
 $ f = f(t), \ t \in [0,1] $ the slight modified module of uniform continuity

$$
\omega(f,\delta) = 0.25 \sup_{ h: |h| \le \delta} \sup_{t \in [0,1]}
| f(t + h) \ - \ f(t)|,
$$

where $ \delta \in [0,1], \ t+h = \min(t+h,1), \ h > 0, $ and
$ t + h = \max(t+h,0) $ in the case $ h < 0. $ \par

 Note that
$$
\lim_{\delta \to 0+} \omega(\xi, \delta) = 0
$$

with probability one and

$$
|\omega(\xi, \delta)| \le 0.5 \ |\xi|_{\infty}.
$$

 It follows from the dominated convergence theorem that

$$
\lim_{\delta \to 0+} {\bf E} \ \Phi( \omega(\xi, \delta)) = 0. \eqno(7)
$$

{\bf 3.} In the language of the theory of Orlicz spaces (see [3], [4], [7] )
 the equality (7) denotes the {\it moment convergence}, or
{\it convergence in mean} the function $ \delta \to \omega(\xi,\delta) $
to zero  as $ \delta \to 0+. $ \par
 If the function $ \Phi $ satisfies the $ \Delta_2 $ condition, the equality
(7) means also  the Orlicz norm convergence, but in general case from
(7) follows only the Orlicz norm $ Or(\Psi) $ convergence:

$$
\lim_{\delta \to 0+} || \omega(\xi, \delta)|| Or(\Psi) = 0. \eqno(8)
$$

{\bf 4.} It follows from (8) then there exists monotonically decreasing
tending to zero sequence  $ \delta(n), \ n = 1,2,\ldots, $ for which

$$
|| \ \omega(\xi, \delta(n) ) \ || Or(\Psi) \le 4^{-n}. \eqno(9)
$$

 Let us define the new subspace  $ Z = (Z, |\cdot|_{\infty, \omega} ) $ of
the space  $ C[0,1] $ consisting on all the (continuous) function $ \{g \} $
for which

$$
\lim_{n \to \infty} 2^n \ \omega(g, \delta(n)) = 0, \eqno(10)
$$

with the finite norm

$$
|g|_{\infty, \omega} = |g|_{\infty} + \sup_n 2^n \  \omega(g, \delta(n)). \eqno(11)
$$

{\bf 5.} Since the functions from the space $ Z $ satisfies the condition (10),
the closed subspace $ Z = (Z, |\cdot|_{\infty, \omega}) $ is SBS. It follows from
the famous Arzela theorem that the space $ Z $ is compact embedded in the space $ C[0,1]. $ \par

{\bf 6.} Finally, let us prove that the Orlicz $ Or(\Psi) $
norm of the  $ \xi $ in the space $ Z $ is finite. We estimate:

$$
|| \ |\xi|_{\infty, \omega} \ ||Or(\Psi) \le || \ |\xi|_{\infty} \ ||Or(\Psi) +
$$

$$
\sum_{n=1}^{\infty} 2^n || \omega(\xi, \delta(n) )|| Or(\Psi) \le
C \ || \ |\xi|_{\infty} \ ||Or(\Phi) + 1 < \infty.
$$

 This completes the proof of theorem 1. \par

\vspace{3mm}

{\bf 3. The case of infinite measure.} \\

\vspace{3mm}

 In this pilcrow we consider the case when the measure $ {\bf P } $ is
unbounded,  but is sigma \ - \ finite and diffuse. \par
 In order to distinct  the probabilistic case and the case with unbounded
measure, we will write the source triple as $ (\Omega, B, {\bf Q}). $\par
 We understood here instead the random variable $ \xi $ some measurable
function  from the triple $ (\Omega, B, {\bf Q} ) $ with Borelian
distribution

$$
\mu_{\xi}(A) = {\bf P}(\xi(\cdot) \in A), \ A \in B.
$$

 Recall that in this case the Orlicz function $ \Psi $ is called
weaker than $ \Phi, \ \Psi << \Phi $ if both the function $ \Psi, \ \Phi $
satisfy the condition (5) and in addition satisfy the following condition:
for all positive constant $ v; \ v = const > 0 $

$$
\lim_{u \to 0}  \frac{\Psi(v \ u ) }{\Phi(u)} = 0.   \eqno(12)
$$

 The next result is proved analogously to the Theorem 1; the important
facts for the moment and norm convergence in Orlicz spaces with unbounded
measure see in the books [3], [4].\par

{\bf Theorem 2.} Let $ \xi $ be measurable function $ \xi: \Omega \to X $
 with Borelian distribution
$ \mu_{\xi}(\cdot) $ in the SBS $ (X, | \cdot|X |) $ satisfies the conditions (1) and (4). \par
 Let also $ \Psi(\cdot) $ be Orlicz function weaker than $ \Phi(\cdot):
\Psi << \Phi. $ \par

 Then there exists a compact embedded in the space $ X $ the SBS space $ Y =
(Y, |\cdot|Y) ), $  which depended on the $ \mu_{\xi} $ and on the function
$ \Psi, $ such that the variable $ \xi $ belongs to the space $ Y $  a.e. and moreover

$$
|| \ |\xi|Y \ ||Or(\Psi) < \infty. \eqno(13)
$$

\vspace{3mm}

{\bf 4. Some generalizations.}\par

\vspace{3mm}

{\bf A.}  The assertion of theorems 1 and 2 is true for the {\it family }
$ M = \{ \mu_{\alpha } \}, \ \alpha \in W, \ W $  is arbitrary set,
of Borelian  $ \sigma \ - $ finite distributions $ \{ \mu_{\alpha } \} $
(i.e. not necessary to be bounded),  if the family $ M $ is dominated in the Radon \ - \ Nikodim sense. \par
 For instance, this is true for countable family $ M. $ \par

 Indeed, there exists a {\it common} compact embedded in the space $ X $
subspace $ \ (Y, \ |\cdot|Y) \ $  such that  for all $ \alpha \in W $

$$
\mu_{\alpha}( X \setminus Y) = 0.
$$

 The assertion of theorem 1 and 2 about Orlicz uniform integrability of a
family $ M $ is also true.\par
 Namely, let for some Orlicz function $ \Phi = \Phi(u) $ 
 
 $$
 \sup_{\alpha \in W} \int_X \Phi(x) \ \mu_{\alpha}(dx) < \infty;
 $$
 and let $ \Psi = \Psi(u) $ be an Orlicz function  weaker than $ \Phi: \ 
 \Psi << \Phi; $ then the common for all values $ \alpha $ compact embedded   subspace $ Y $  may be constructed such that 

 $$
 \sup_{\alpha \in W} \int_Y \Psi(y) \ \mu_{\alpha}(dy) < \infty.
 $$

{\bf B.} Let $ M = \{ \mu_m \}, \ m = 1,2,\ldots $ be countable family of
Borelian distributions in SBS $ \ (X, \ |\cdot|X) \ $ which convergent weakly
to  some distribution $ \mu, \ $ i.e. for all continuous bounded functional 
$ g: X \to R $

$$
\lim_{n \to \infty} \int_X g(x) \ \mu_n(dx) = \int_X g(x) \ \mu(dx).
$$

  For example, the family $ M $ may satisfy the
 to \ - \ called Central Limit Theorem  (CLT) in the space $ X. $ \par
 This means by definition that the measure $ \ \mu \ $ is Gaussian.\par

 Then the common compact embedded  support subspace $ \ (Y, \ |\cdot|Y) ) \ $ can be constructed in addition such that the sequence $ M $ convergent weakly 
to the measure $ \mu $ in the space $ Y: $  for all continuous bounded 
functional $ h: Y \to R $

$$
\lim_{n \to \infty} \int_Y h(y) \ \mu_n(dy) = \int_Y h(y) \ \mu(dy).
$$

{\bf C.} V.V.Buldygin [2] proved that in the probabilistic case $ \mu(X) = 1 $
the subspace $ Y $ may be constructed to be reflexive and with continuous
differentiable in the Freshe sense norm. \par
 At the same is true also in general, i.e. in unbounded case $ \mu(X) =
\infty.$ \par

{\bf D.} For the space $ X = C(T), $ where $ T $ is compact metric space
with distance $ d = d(t,s), \ t,s \in T, $
the assertion of Theorem 1 may be formulated as follows. \par
 Let $ \xi = \xi(t), \ t \in T $ be continuous with probability one random
field such that  for some Orlicz function $ \Phi $

$$
{\bf E } \Phi \left(|\xi|_{\infty} \right) < \infty.
$$
 Let also $ \Psi $ be another Orlicz function, $ \Psi << \Phi. $ \par
 Then there exist a  non \ - \ negative
r.v. $ \zeta $ and non \ - \ random continuous semi \ - \
distance $ r = r(t,s) $ on the space $ T $ such that

$$
|\xi(t) - \xi(s)| \le \zeta \times r(t,s),
$$

where

$$
{\bf E} \Psi(\zeta) < \infty.
$$

{\bf DH.} In the connection of the last assertion we dare formulate
the following, interest by our opinion,  hypothesis. Let $ \theta =
\theta(t), \ t \in T $ be arbitrary separable random field, 
 centered: $ \ {\bf E} \theta(t) = 0 \ $ or not, bounded with probability 
one:

$$
\sup_{t \in T} |\theta(t) | < \infty \ {\bf a.e.}.
$$

 Assume in addition that for some Orlicz function $ \Phi(\cdot) $

$$
\sup_{t \in T} ||\theta(t)||Or(\Phi) < \infty.
$$

 Let also $ \Psi $ be another Orlicz function such that $ \Psi << \Phi. $ \par

 Open question: there holds (or not)

$$
|| \ \sup_{t \in T} |\theta(t)| \ ||Or(\Psi) < \infty ? \eqno(14)
$$

 The conclusion (14) is true for the {\it centered} Gaussian fields [8]; if the field
$ \theta(\cdot) $ satisfies the so \ - \ called  entropy or generic chaining conditions [5], [9], [10]; in the case if $ \theta $ belongs to the
domain of attraction of Law of Iterated Logarithm [11] etc. \par
  Finally, let us consider the following example.
  Let $ \tau $ be Normal (Gaussian) standard distributed r.v.: $ Law(\tau) = N(0,1) $
 ant let $ T = R = (-\infty, \infty). \ $ We define

$$
\theta(t) = \tau \cdot t \ -  \ |t|^p/p, \ p = const \in (1, 2).
$$
  Then $ \theta(\cdot) $ is {\it upper \ - \  bounded  } Gaussian {\it non \ - \
  centered } random process, but

$$
\sup_{t \in R} \theta(t) = |\tau|^q/q, \ q = p/(p \ - \ 1) > 2.
$$

 The tail of distribution of the r.v. $ \sup_{t \in R} \ \theta(t)  $ is essentially
heaver in comparison to the Gaussian r.v. $ \tau $  or to the upper tail of each variable
$ \theta(t), \ t \in T. $ \\

\vspace{3mm}

\begin{center}

{\bf REFERENCES }\\

\end{center}

  1.  {\sc Ostrovsky E.I.} (1980). {\it On the support of probabilistic measures in separable Banach spaces. } Soviet Mathematic,  Doklady, v.255, $ N^o \ 6 $ pp. 836
  \ - \ 838, \ (in Russian).\\

   2. {\sc Buldygin V.V.} (1984). {\it  Supports of probabilistic measures in
   separable Banach spaces.} Theory Probab. Appl. {\bf 29} v.3, pp. 528 \ - 532,
   \ (in Russian). \\

   3. {\sc  M.M. Rao, Z.D.Ren.} (1991). {\it Theory of Orlicz Spaces.}  Basel - New   York, Marcel Decker. \\

   4. {\sc M.M. Rao, Z.D.Ren.} (2002). {\it Application of Orlicz Spaces.}  Basel - New York, Marcel Decker. \\

   5. {\sc Ostrovsky E.I.} (1999). {\it Exponential estimations for Random
    fields and its applications } (in Russian). Moscow - Obninsk, OINPE.\\

  6. {\sc Ostrovsky E.I.} (2002). {\it Exact exponential estimations for
     random field maximum distribution.} Theory Probab. Appl. {\bf 45} v.3,
      pp. 281 - 286, (in Russian).  \\

  7. {\sc M.A.Krasnoselsky, Ya.B.Rutisky.} (1961). {\it Convex functions and Orlicz's Spaces.} P. Noordhoff LTD, The Netherland, Groningen. \\

  8. {\sc Fernique X.} (1975). {\it Regularite des trajectoires des
    function aleatiores gaussiennes.} Ecole de Probablite de
    Saint-Flour, IV – 1974, Lecture Notes in Mathematic. {\bf 480} 1 – 96, Springer
    Verlag, Berlin.\\

   9. {\sc Talagrand M.} (1996). {\it Majorizing measure: The generic chaining.}
       Ann. Probab. {\bf 24} 1049 - 1103. MR1825156 \\

   10. {\sc Ostrovsky E., Rogover E.} {\it Exact exponential bounds for the random
    Fields Maximum Distribution via the majoring Measures (generic Chaining)}.
     Electronic Publications, arXiv:o802v1 [math.PR], 4 Feb 2008.\\

   11. {\sc Ostrovsky E.I.} {\it Exponential Bounds  in the Law of Iterated
     Logarithm in Banach Space.} (1994), Math. Notes, {\bf 56}, 5, p. 98 - 107. \\

\end{document}